# The Magnetic Tower of Hanoi

*Uri Levy*
*Atlantium Technologies, Har-Tuv Industrial Park, Israel*

uril@atlantium.com

**Abstract**

In this work I study a modified Tower of Hanoi puzzle, which I term Magnetic Tower of Hanoi (MToH). The original Tower of Hanoi puzzle, invented by the French mathematician Edouard Lucas in 1883, spans "base 2". That is – the number of moves of disk number k is $2^{(k-1)}$, and the total number of moves required to solve the puzzle with N disks is $2^N - 1$. In the MToH puzzle, each disk has two distinct-color sides, and disks must be flipped and placed so that no sides of the same color meet. I show here that the MToH puzzle spans "base 3" - the number of moves required to solve an N+1 disk puzzle is essentially three times larger than he number of moves required to solve an N disk puzzle. The MToH comes in 3 flavors which differ in the rules for placing a disk on a free post and therefore differ in the possible evolutions of the Tower states towards a puzzle solution. I analyze here algorithms for minimizing the number of steps required to solve the MToH puzzle in its different versions. Thus, while the colorful Magnetic Tower of Hanoi puzzle is rather challenging, its inherent freedom nurtures mathematics with remarkable elegance.

**The Classical Tower of Hanoi**

The classical Tower of Hanoi (**ToH**) puzzle[1,2,3] consists of three posts, and N disks. The puzzle solution process ("game") calls for one-by-one disk moves restricted by one "size rule". The puzzle is solved when all disks are transferred from a "Source" Post to a "Destination" Post.

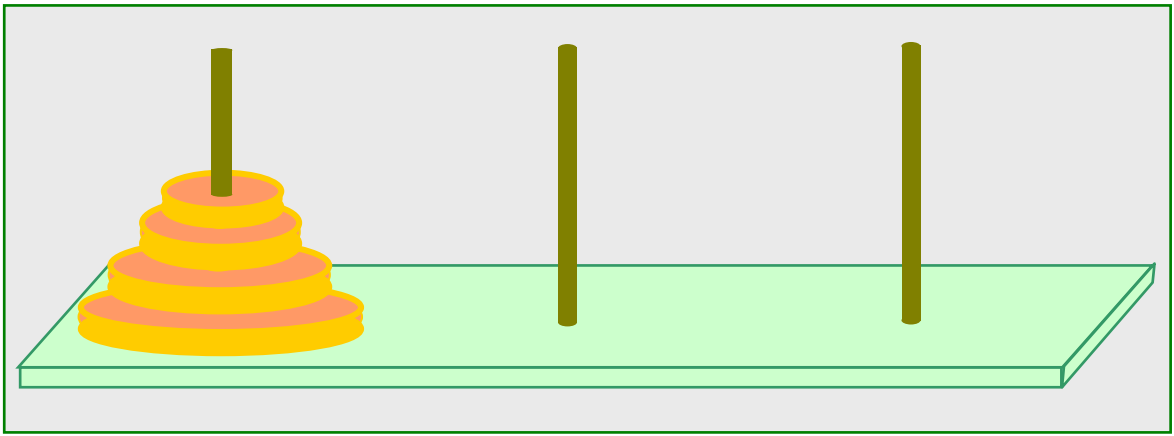

**Figure 1:** *The classical Tower of Hanoi puzzle. The puzzle consists of three posts, and N disks. The puzzle solution process ("game") calls for one-by-one disk moves restricted by one "size rule". The puzzle is solved when all disks are transferred from a "Source" Post to a "Destination" Post. The minimum number of disk-moves necessary to solve the **ToH** puzzle with N disks is $2^N - 1$.*

Let's define the **ToH** puzzle in a more rigorous way.

### 1.1. The Classical Tower of Hanoi – puzzle description

A more rigorous description of the **ToH** puzzle is as follows -

**Puzzle Components:**
- Three equal posts
- A set of N different-diameter disks

**Puzzle-start setting:**
- N disks arranged in a bottom-to-top descending-size order on a "Source" Post (Figure 1)

**Move:**
- Lift a disk off one Post and land it on another Post

**Disk-placement rules:**
- **The Size Rule:** A small disk can not "carry" a larger one (Never land a large disk on a smaller one)

**Puzzle-end state:**
- N disks arranged in a bottom-to-top descending-size order on a "Destination" Post (one of the two originally-free posts)

Given the above description of the classical **ToH** puzzle, let's calculate the (minimum) number of moves necessary to solve the puzzle.

### 1.2. Number of moves

Studying the classical **ToH** puzzle in terms of required moves to solve the puzzle, it is not too difficult to show[2,3] (prove) that the k-th disk will make $P(k)$ moves given by

$$P(k) = 2^{k-1}. \tag{1}$$

Disk numbering is done of course from bottom to top and as Equation 1 states - the smallest disk ($k = N$), like the least significant bit in a "binary speedometer", makes the largest number of moves.

The total number of moves $S(N)$ is given by the sum -

$$S(N) = \sum_{k=1}^{N} 2^{k-1} = 2^N - 1. \tag{2}$$

Table 1 lists the (minimum) number of moves of each disk (Equation 1) and the total (minimum) number of moves required to solve the classical **ToH** puzzle (Equation 2) for the first eight stack heights.

| N \ k | 1 | 2 | 3 | 4 | 5 | 6 | 7 | 8 | SUM | $2^N - 1$ |
|---|---|---|---|---|---|---|---|---|---|---|
| 1 | 1 | | | | | | | | 1 | 1 |
| 2 | 1 | 2 | | | | | | | 3 | 3 |
| 3 | 1 | 2 | 4 | | | | | | 7 | 7 |
| 4 | 1 | 2 | 4 | 8 | | | | | 15 | 15 |
| 5 | 1 | 2 | 4 | 8 | 16 | | | | 31 | 31 |
| 6 | 1 | 2 | 4 | 8 | 16 | 32 | | | 63 | 63 |
| 7 | 1 | 2 | 4 | 8 | 16 | 32 | 64 | | 127 | 127 |
| 8 | 1 | 2 | 4 | 8 | 16 | 32 | 64 | 128 | 255 | 255 |

**Table 1:** *Minimum number of disk-moves required to solve the classical Tower of Hanoi puzzle. N is the total number of disks participating in the game and k is the disk number in the ordered stack, counting from bottom to top. The k-th disk "makes" $2^{(k-1)}$ moves (Equation 1). The total number of disk-moves required to solve an N-disk puzzle is $2^N - 1$ (Equation 2).*

Table 1 clealy shows how (elegantly) the classical **ToH** spans base 2.

Let's see now how base 3 is spanned by the far more intricate Magnetic Tower of Hanoi puzzle.

## 2. The Magnetic Tower of Hanoi

In the Magnetic Tower of Hanoi puzzle[4], we still use three posts and N disks. However, the disk itself, the move definition and the game rules are all modified (extended).

The rigorous description of the **MToH** puzzle is as follows:

**Puzzle Components:**
- Three equal posts
- A set of N different-diameter disks
- Each disk's "bottom" surface is colored **Blue** and its "top" surface is colored **Red**

**Puzzle-start setting:**
- N disks arranged in a bottom-to-top descending-size order on a "Source" Post (Figure 2)

- The **Red** surface of every disk in the stack is facing upwards (Figure 2). *Note that the puzzle-start setting satisfies the "Magnet Rule" (see below). And needless to say, **Red** is chosen arbitrarily without limiting the generality of the discussion.*

**Move:**
- Lift a disk off one post
- Turn the disk upside down and land it on another post

**Disk-placement rules:**
- **The Size Rule:** A small disk can not "carry" a larger one (Never land a large disk on a smaller one)
- **The Magnet Rule:** Rejection occurs between two equal colors (Never land a disk such that its bottom surface will touch a co-colored top surface of the "resident" disk)

**Puzzle-end state:**
- N disks arranged in a bottom-to-top descending-size order on a "Destination" Post (one of the two originally-free posts)

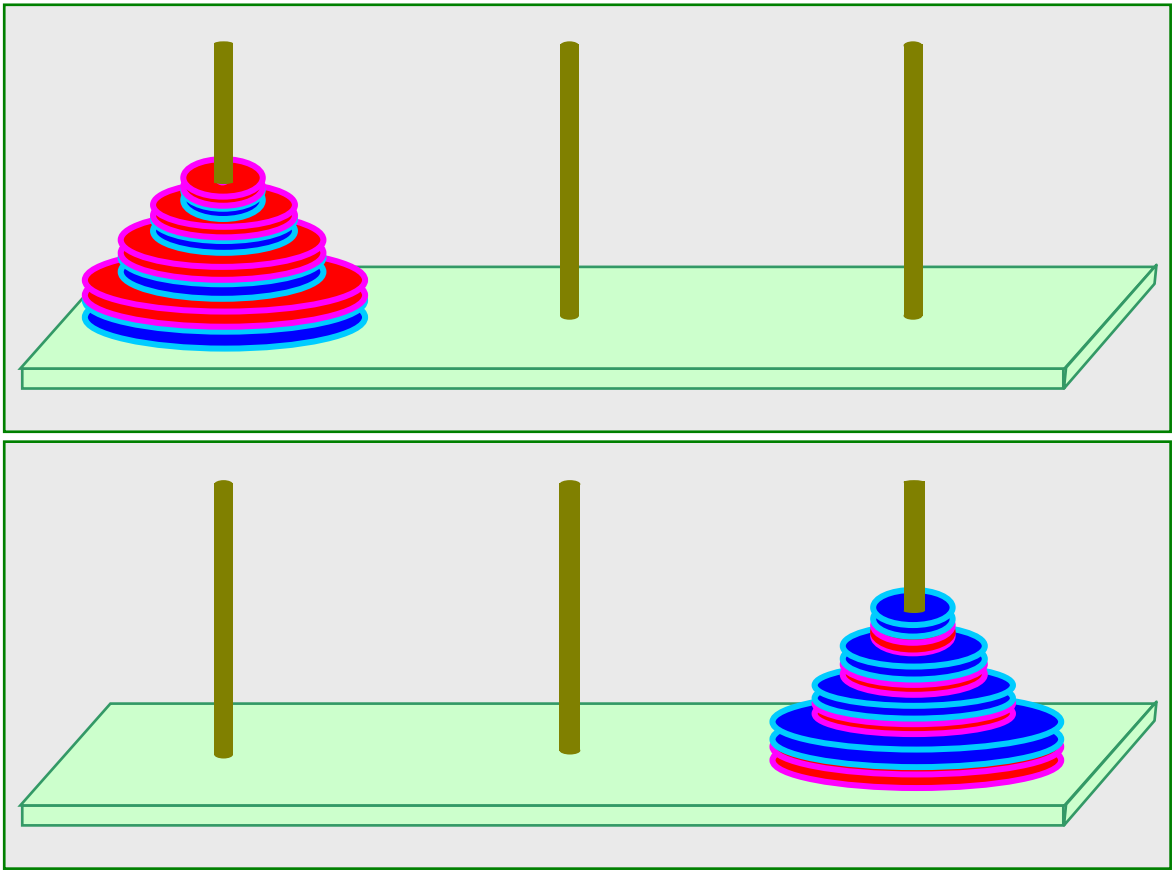

**Figure 2:** *The Magnetic Tower of Hanoi puzzle. Top – puzzle-start setting. The puzzle consists of three posts, and N two-color disks. The puzzle solution process ("game") calls for one-by-one disk moves restricted by two rules – the Size Rule and the Magnet Rule. The puzzle is solved when all disks are transferred from a "Source" Post to a "Destination" Post - bottom.*

Given the above description of the **MToH** puzzle, let's calculate the number of moves necessary to solve the puzzle.

We start by explicitly solving the N=1, N=2 and N=3 cases.

## 2.1. Explicit solution for the first three stacks of the MToH puzzle

**The N = 1 case** is trivial – move the disk from the Source Post to a Destination Post (Figure 3).

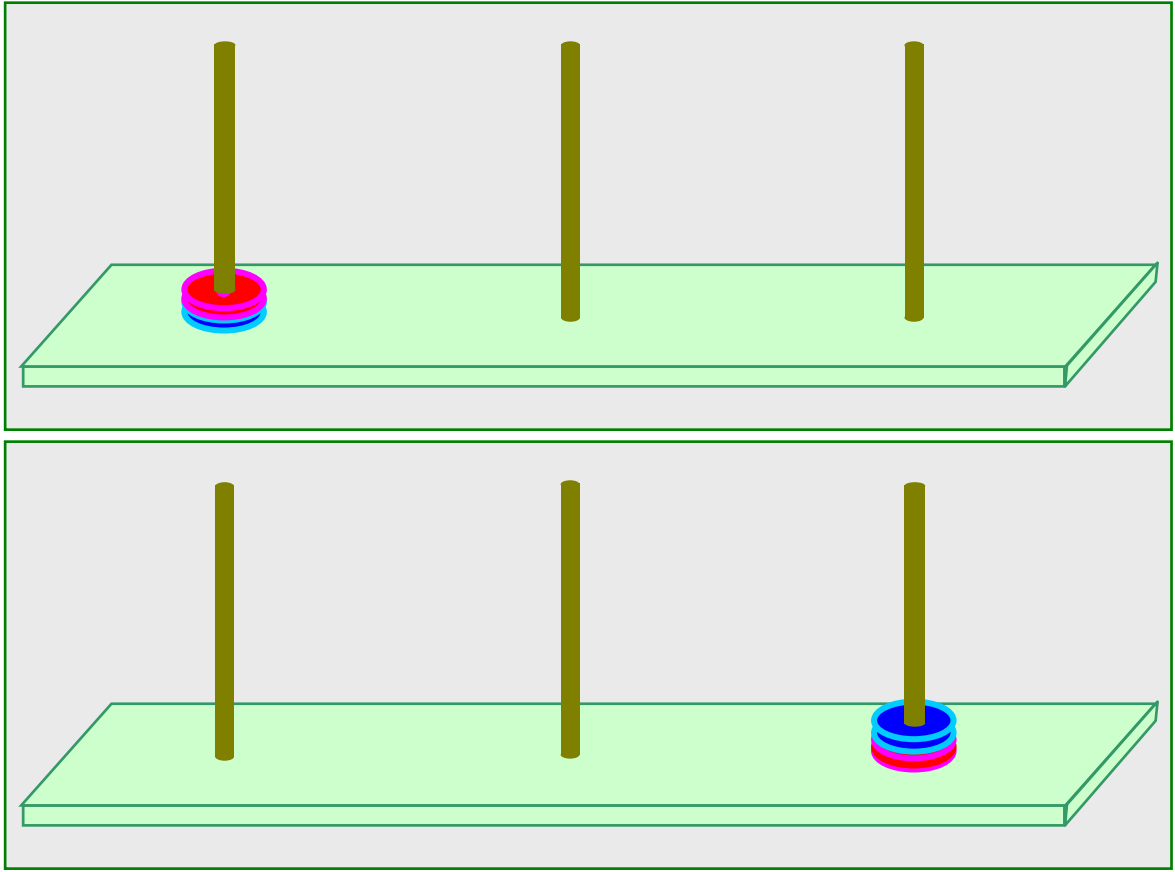

**Figure 3:** *The start-setting (top) and the end-state (bottom) for the N=1 MToH puzzle. The number of moves required to solve the puzzle is P(1) = 1.*

Thus, for the N=1 case we have

$$P(1) = 1 \; ; \; S(1) = 1. \tag{3}$$

Let's see **the N=2 case**.

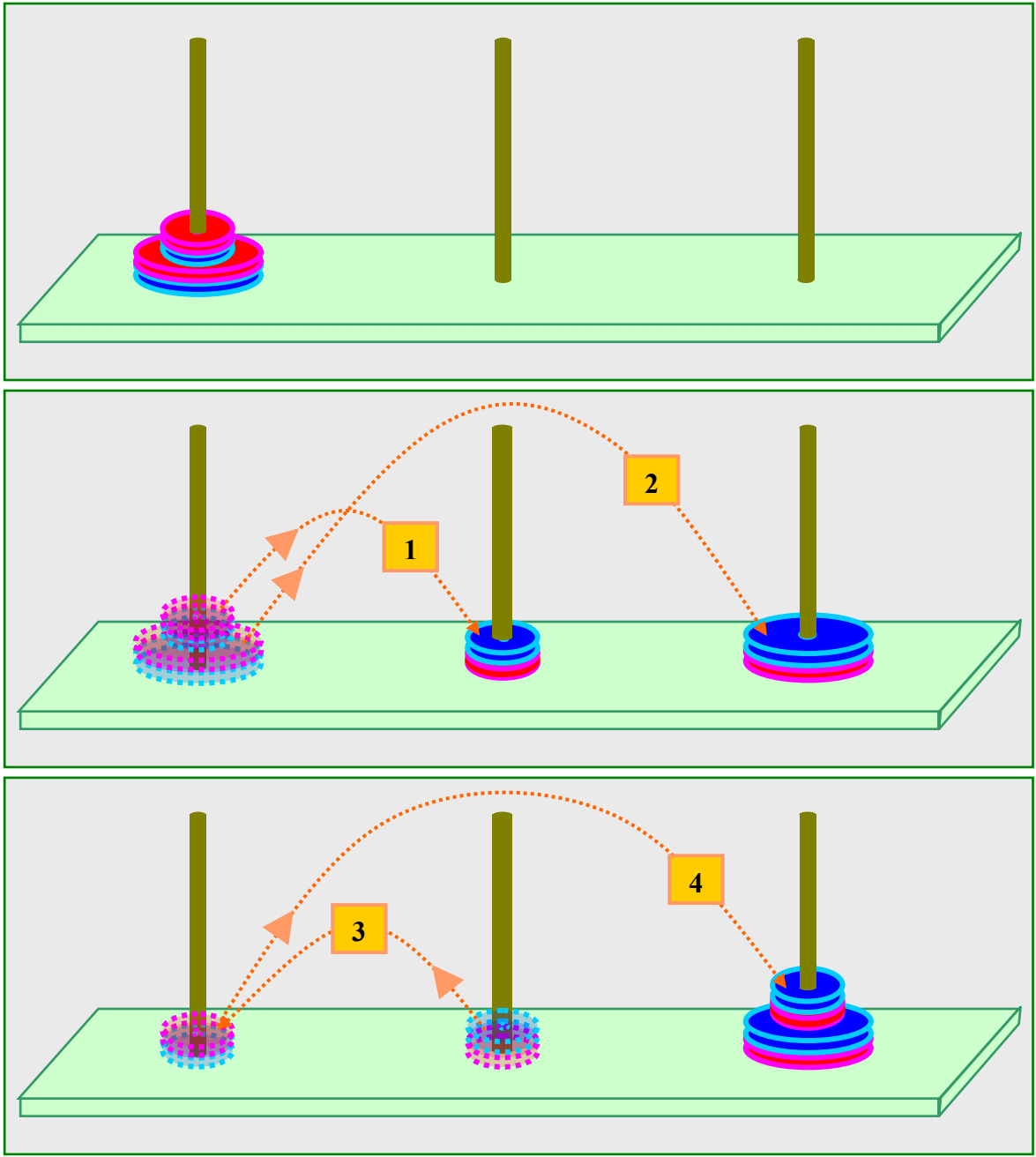

**Figure 4:** *The start-setting (top), an intermediate setting (center) and the end-state (bottom) for the N=2 **MToH** puzzle. The number of moves to progress from the start-setting to the intermediate state described by the center figure is 2. The number of moves to progress from the center-described state to the end-state described by the bottom figure is again 2. Thus, the (minimum) number of moves required to solve the puzzle is S(2) = 4. Note that two different solution routes, both of length 4, exist (1,2,1,1 – shown, 1,1,2,1 – not shown).*

Consulting Figure 4 we find for the N=2 case -

$$S(2) = 4 \, .\tag{4}$$

The small disk made 3 (=$3^1$) moves and the large disk made 1 (=$3^0$) move. Thus far then, for the N=1 and N=2 cases, base 3 is elegantly spanned as

$$P(k) = 3^{k-1} \text{ and } S(N) = (3^N - 1)/2 \text{ ; N = 1,2.}\tag{5}$$

Exactly analogous to the base 2 span by the classical **ToH** (Equation 1 and Equation 2).

But let's see now **the N=3 case**.

To conveniently talk about the N=3 case, let's (arbitrarily and without loss of generality) define the posts (refer to Figure 5 below) as

- S – the Source Post
- D – the Destination Post
- I – the (remaining) Intermediate Post

Let's also memorize the disk numbering convention:

- 1 – the largest disk
- 2 – the mid-size disk
- 3 – the smallest disk

| Step number | Disks | From | To | # of moves | Comments |
|---|---|---|---|---|---|
| 1 | 2,3 | S (Red) | I (Blue) | 4 | Equation 4 |
| 2 | 1 | S (Red) | D (Blue) | 1 | S now free |
| 3 | 3 | I (Blue) | D (Blue) | 2 | Through S, S now free |
| 4 | 2 | I (Blue) | S (Red) | 1 | Post I now free<br>Fig. 5 - middle |
| 5 | 3 | D (Blue) | I (Red) | 1 | S and I are both Red<br>I switched Blue→Red |
| 6 | 2 | S (Red) | D (Blue) | 1 |  |
| 7 | 3 | I (Red) | D (Blue) | 1 | Puzzle solved |
|  |  |  |  | 11 | **Total # of moves** |

**Table 2:** *Explicit description of the moves to solve the N=3 **MToH** puzzle. The total number of moves is 11, which does NOT exactly coincide with the "base 3 rule" (Equation 5).*

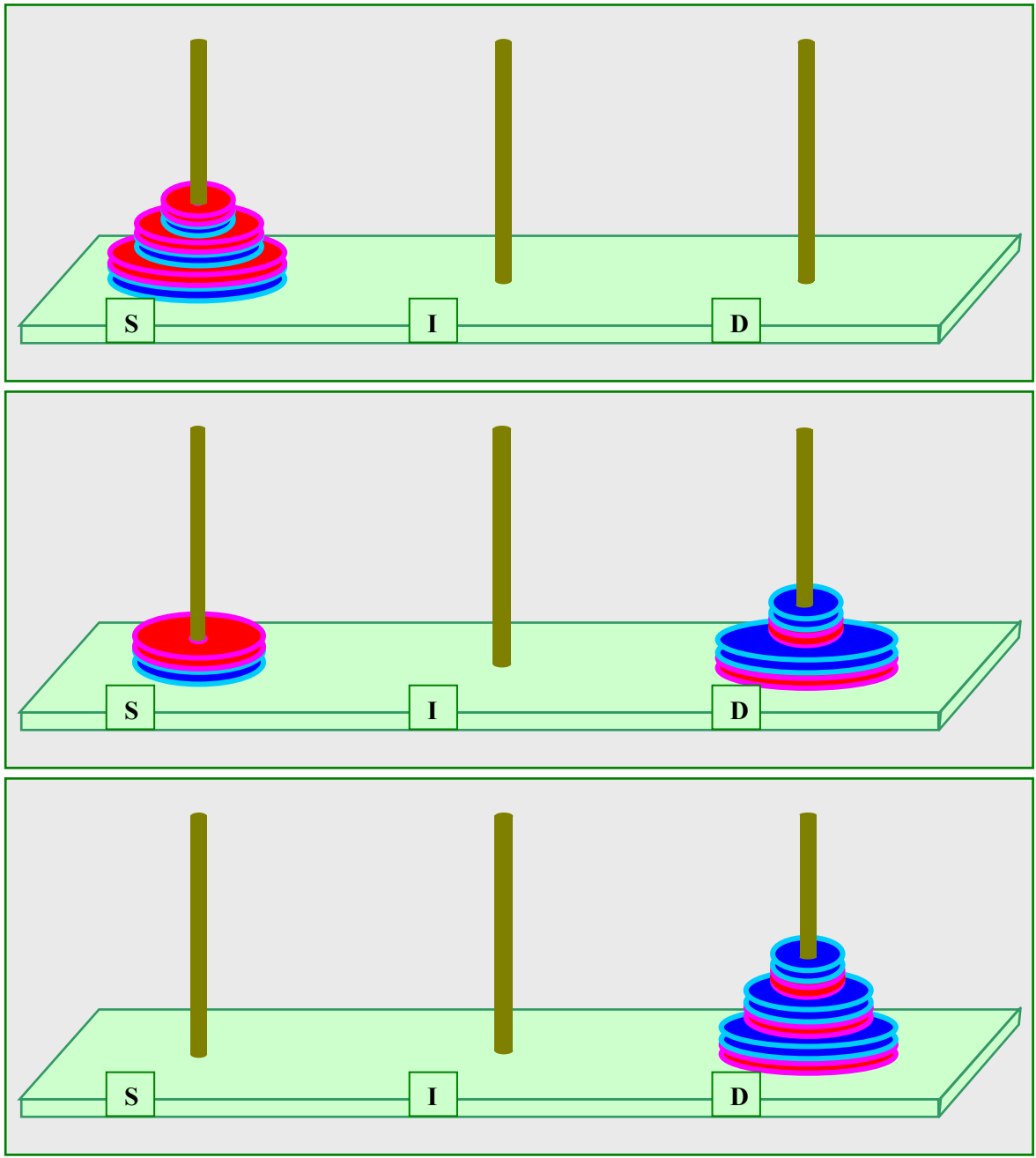

**Figure 5:** *The start-setting (top), an intermediate setting (center) and the end-state (bottom) for the N=3 **MToH** puzzle. The number of moves to progress from the start-setting to the intermediate state described by the center figure is 8 (read the text for details). The number of moves to progress from the center-described state to the end-state described by the bottom figure is 3. Thus, the number of moves required to solve the puzzle is S(3) = 11. The S(3) number (11) breaks the perfect base 3 rule (Equation 5). We therefore need to probe into the puzzle further in order to decipher the mystery of this newly observed irregularity (and come up with a modified rule).*

Listed in table 2 are the moves required to solve the N=3 **MToH** puzzle.

As shown in Table 2 and as demonstrated by Figure 5 –

$$S(3) = 11. \tag{6}$$

The resulting total number of moves violates the "base 3 rule" (should have been 13, refer to Equation 5). The states of the puzzle before step 1 (puzzle-start state), after step 4 and after step 7 (puzzle-end state) are shown by Figure 5 – top, center, bottom respectively.

In order to decipher the mystery (of the newly observed irregularity) and to progress with the analysis of the **MToH** puzzle, let's define a new magnetic tower, refer to it as the "Colored Magnetic Tower of Hanoi" and study its properties.

### 2.2. The Colored Magnetic Tower of Hanoi – the "100" solution

Studying the N=3 **MToH** puzzle, I realized that what breaks the base 3 rule is the possibility of the smallest disk to move to a free post (step 5 in Table 2). By "free" I mean a post that is not "magnetized" or not "color coded". A Neutral Post that can accept any-color disk. To suppress this freedom, let's permanently color-code each post, call the restricted tower the Colored Magnetic Tower of Hanoi (**CMToH**) and see what happens.

#### 2.2.1. Definition of "Colored"

Let's start with a definition of a **C**olored **MToH**:

An **MToH** is "Colored" if (without loss of generality) its posts are pre-colored (or "permanently colored")

**Either as**

        1. Red-Blue-Blue

**Or as**

        2. Red-Red-Blue.

Let's designate this (permanently) Colored **MToH as CMToH**.

The two versions of the newly defined **CMToH** puzzle are shown by Figure 6. The moves to solve the **CMToH** puzzle with N=2, for each of its versions, are explicitly detailed by Table 3. Note that the only difference between the versions is the "timing" of the move of the big disk (after one move of the small disk in the first version and after two moves of the small disk in the second version).

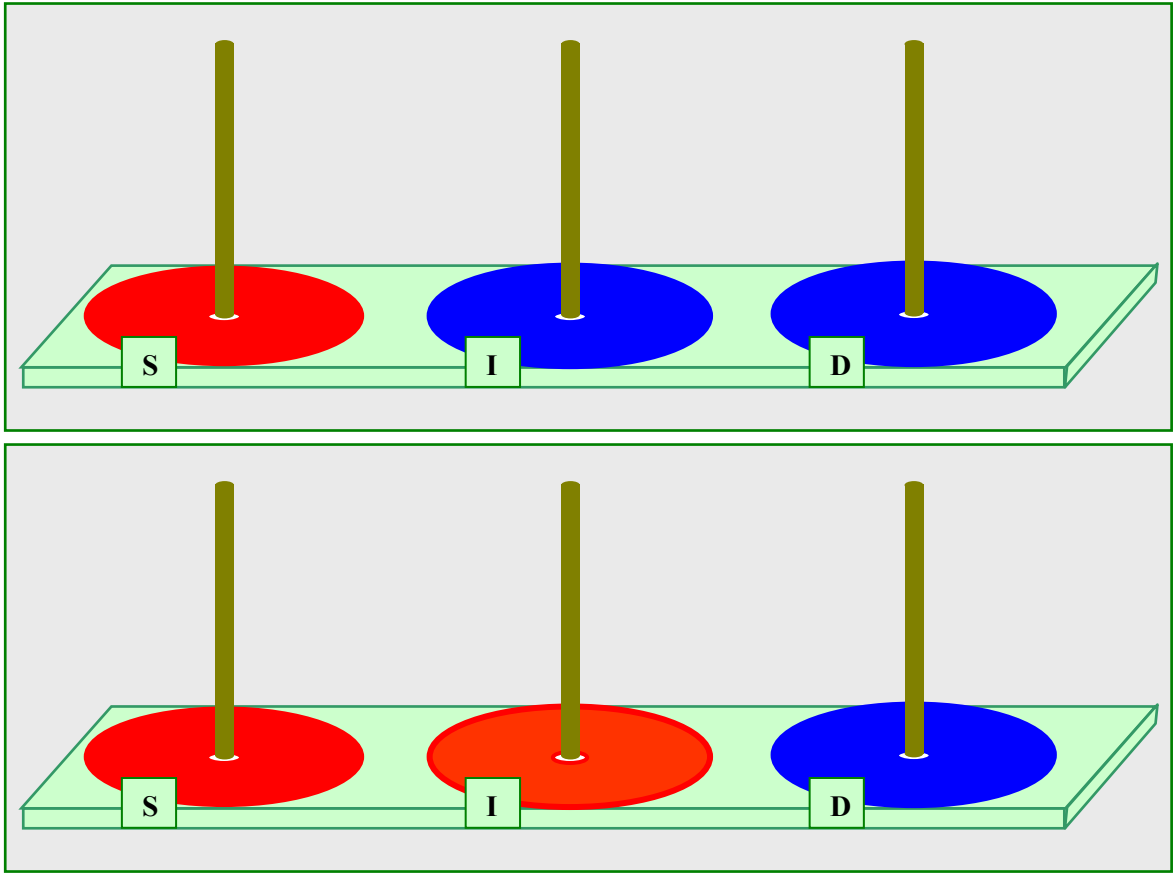

**Figure 6:** *The two versions of the (permanently) Colored Magnetic Tower of Hanoi. As shown in the text, the two definitions are equivalent in terms of number of moves. Given a Colored Magnetic Tower of Hanoi, the number of moves of disk k are $P(k) = 3^{(k-1)}$ and the total number of moves is $S(N) = (3^N – 1)/2$. Thus, the freshly defined **C**olored **M**agnetic **T**ower **of H**anoi strictly spans base 3.*

| Step number | Disks | From | To | # of moves | Comments |
|---|---|---|---|---|---|
| | | | | | **1. Red-Blue-Blue** |
| 1 | 2 | S (Red) | I (Blue) | 1 | |
| 2 | 1 | S (Red) | D (Blue) | 1 | |
| 3 | 2 | I (Blue) | D (Blue) | 2 | Through **Red** S |
| | | | | | **2. Red-Blue-Red** |
| 1 | 2 | S (Red) | I (Red) | 2 | Through **Blue** D |
| 2 | 1 | S (Red) | D (Blue) | 1 | |
| 3 | 2 | I (Red) | D (Blue) | 1 | |

**Table 3:** *Explicit description of the moves to solve the N=2 **CMToH** puzzle. The total number of moves for both versions is 4. And the N=3 case of the **CMToH** puzzle is solved by 13 moves.*

### 2.2.2. Expressions for the number of moves

Simple observations reveal that, as is the case with the classical **ToH**, the **"forward" moves solving the CMToH puzzle are deterministic**.

Furthermore, it is not too difficult to show by a recursive argument (see the proof for the classical **ToH**[2,3]) that the number of disk moves $P_{100}(k)$ and (therefore) the total number of moves $P_{100}(N)$ perfectly span base 3:

$$P_{100}(k) = 3^{k-1}. \qquad (7)$$

and

$$S_{100}(N) = \sum_{k=1}^{N} 3^{k-1} = \frac{3^N - 1}{3^1 - 1}. \qquad (8)$$

The subscript "100" in Equations 7 and 8, relates to a solution "duration" of 100%.

Table 4 lists the (minimum) number of moves of each disk (Equation 7) and the total (minimum) number of moves (Equation 8) required to solve the **CMToH** puzzle for the first eight stack heights.

| N \ k | 1 | 2 | 3 | 4 | 5 | 6 | 7 | 8 | SUM | $(3^N - 1)/2$ |
|---|---|---|---|---|---|---|---|---|---|---|
| 1 | 1 | | | | | | | | 1 | 1 |
| 2 | 1 | 3 | | | | | | | 4 | 4 |
| 3 | 1 | 3 | 9 | | | | | | 13 | 13 |
| 4 | 1 | 3 | 9 | 27 | | | | | 40 | 40 |
| 5 | 1 | 3 | 9 | 27 | 81 | | | | 121 | 121 |
| 6 | 1 | 3 | 9 | 27 | 81 | 243 | | | 364 | 364 |
| 7 | 1 | 3 | 9 | 27 | 81 | 243 | 729 | | 1093 | 1093 |
| 8 | 1 | 3 | 9 | 27 | 81 | 243 | 729 | 2187 | 3280 | 3280 |

**Table 4:** *Minimum number of disk-moves required to solve the Colored Magnetic Tower of Hanoi puzzle. N is the total number of disks participating in the game and k is the disk number in the ordered stack, counting from bottom to top. The k-th disk "makes" $3^{(k-1)}$ moves (Equation 7). The total number of disk-moves required to solve an N-disk puzzle is $(3^N – 1)/2$ (Equation 8).*

Having solved the rather simple Colored Magnetic Tower puzzle, we can move on to solving the more intricate "Free" or "Dynamically Colored" Magnetic Tower puzzle. As discussed below, we will identify "Free" and "Colored" states of the "Dynamically Colored" **MToH** leading to a far "shorter" solution (relative to the "100" solution) of the **MToH** puzzle.

### 2.3. The "67%" solution of the MToH puzzle

The color of posts in the **MToH** puzzle is determined by the color of the disks it holds. The color is therefore "dynamic". During the game, the color of a given post can be **RED**, can be **BLUE**, and can be **Neutral**. For moves analysis, we can distinguish between three distinct **MToH** states.

#### 2.3.1. Distinct states of the Magnetic Tower of Hanoi

After playing with the **MToH** puzzle for a while, one may realize that actually three distinct tower states exist

- "Free" - two posts are Neutral ("start" and "end" states)
- "Semi-Free" – one post Neutral, the other two are oppositely Colored
- "Colored" – two posts are co-colored

During the "game", the tower is some-times "Semi-Free". Which opens the room for significant "savings".

The "67" solution indeed takes advantage of the Tower's occasional "semi-freedom". In fact, as I will show below, for large N, the number of moves for this "67" solution is 2/3 of the number of moves of the "100" solution.

#### 2.3.2. The 67% solution

The "67" (percent) solution is based on the sequence listed in Table 5:

| Step # | Disks | From | To | # of moves | Comments |
|---|---|---|---|---|---|
| 1 | 2 to N | S (Red) | I (Blue) | $S_{67}(N-1)$ | "Free" **MToH** |
| 2 | 1 | S (Red) | D (Blue) | 1 | |
| 3 | 3 to N | I (Blue) | D (Blue) | $2*S_{100}(N-2)$ | Through **"Red"** S |
| 4 | 2 | I (Blue) | S (Red) | 1 | |
| 5 | 3 to N | D (Blue) | I (RED) | $1*S_{100}(N-2)$ | |
| 6 | 2 | S (Red) | D (Blue) | 1 | |
| 7 | 3 to N | I (RED) | D (Blue) | $1*S_{100}(N-2)$ | Puzzle solved |

**Table 5:** *The sequence of moves for the "67" solution of the MToH puzzle. As shown, the total number of moves is*
$$S_{67}(N) = S_{67}(N-1) + 4*S_{100}(N-2) + 3.$$

If we start on a Red post (up-facing surfaces of all disks are colored Red), then after $S_{67}(N-1) + 1$ moves we arrive at the state described by Figure 7. The rest of the move-sequence to solve the puzzle is $4*S_{100}(N-2) + 2$, as detailed in the text and as listed in Table 5.

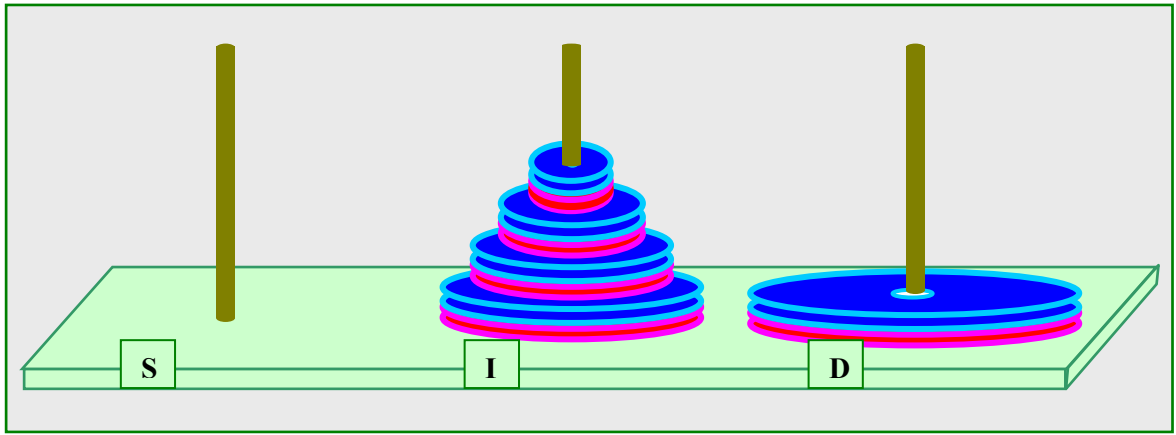

**Figure 7:** *Moving N-disks from S to D by the "67" Algorithm. The figure shows state of the tower (started Red on Post S) after N-1 disks are moved to the Intermediate Post and the N-th disk is moved to the Destination Post. The number of (minimum) moves to get from puzzle-start state to the figure-described state is $S_{67}(N-1) + 1$. The rest of the move-sequence to solve the puzzle is $4*S_{100}(N-2) + 2$, as detailed in the text and as listed in Table 5.*

The recursive proof of the "67" Algorithm is the following: we know how to solve for 3 disks. For $N > 3$, if the algorithm works for N disks, it works for N+1 disks because after we have successfully moved N disks ("down") from S to I (as assumed) and moved the N+1 disk from S to D in a legal way (Figure 7), we move N-2 disks using the always legal "Colored" algorithm (steps 3, 5, 7 in Table 5) and move the N-1 single disk twice in a legal way (steps 4 and 6 in Table 5).

The number of moves in the "67" solution Algorithm as explained above is

$$S_{67}(N) = S_{67}(N-1) + 4 \cdot S_{100}(N-2) + 3 .\qquad(9)$$

Given Equation 7 and Equation 9, we can quickly formulate non-recursive expressions to the number of moves. Consulting these two equations and performing some algebraic manipulations, we find for the "67" solution of the Magnetic Tower of Hanoi –

$$P_{67}(1) = 1$$

$$P_{67}(k) = 2 \cdot 3^{k-2} + 1 \quad ; \quad k \geq 2 \tag{10}$$

And from Equation 10:

$$S_{67}(N) = 1 + \sum_{k=2}^{N} \left[ 2 \cdot 3^{k-2} + 1 \right] = 3^{N-1} + N - 1. \tag{11}$$

Just copying Equation 11 for clarity:

$$S_{67}(N) = 3^{N-1} + N - 1. \tag{12}$$

Table 6 lists the number of moves of each disk (Equation 10) and the total number of moves (Equation 12) for the "67" solution of the **MToH** puzzle, for the first eight stack heights.

| N \ k | 1 | 2 | 3 | 4 | 5 | 6 | 7 | 8 | SUM | $3^{(N-1)} + N-1$ |
|---|---|---|---|---|---|---|---|---|---|---|
| 1 | 1 | | | | | | | | 1 | 1 |
| 2 | 1 | 3 | | | | | | | 4 | 4 |
| 3 | 1 | 3 | 7 | | | | | | 11 | 11 |
| 4 | 1 | 3 | 7 | 19 | | | | | 30 | 30 |
| 5 | 1 | 3 | 7 | 19 | 55 | | | | 85 | 85 |
| 6 | 1 | 3 | 7 | 19 | 55 | 163 | | | 248 | 248 |
| 7 | 1 | 3 | 7 | 19 | 55 | 163 | 487 | | 735 | 735 |
| 8 | 1 | 3 | 7 | 19 | 55 | 163 | 487 | 1459 | 2194 | 2194 |

**Table 6:** *Number of disk-moves for the "67" solution of the "Free" Magnetic Tower of Hanoi puzzle. N is the total number of disks participating in the game and k is the disk number in the ordered stack, counting from bottom to top. The k-th disk "makes" [2\*3$^{(k-2)}$+1] moves (Equation 10). The total number of disk-moves in the "67" solution of the **MToH** puzzle is [3$^{(N-1)}$+N-1] (Equation 12).*

With Equation 8 for the number of moves in the "100" solution and Equation 12 for the number of moves in the "67" solution, one can easily determine the limit of the "duration-ratio" for large stacks:

$$\frac{S_{67}(N)}{S_{100}(N)} = \frac{2}{3} \cdot \frac{1+(N-1)/3^{(N-1)}}{1-1/3^N} \xrightarrow{n \to \infty} \frac{2}{3}. \tag{13}$$

So for large stacks of disks, the duration of the "67" solution is indeed 67% of the duration of the "100" solution.

Knowing the expressions for the exact number of moves for the "100" solution as well as for the "67" solution, we plot a "duration-ratio" curve or "efficiency" curve for the "67" solution – Figure 8. As shown, the curve monotonically (and "quickly") approaches its limit of 2/3 (Equation 13) and with a stack of only seven disks the efficiency curve is practically at its large-number limit.

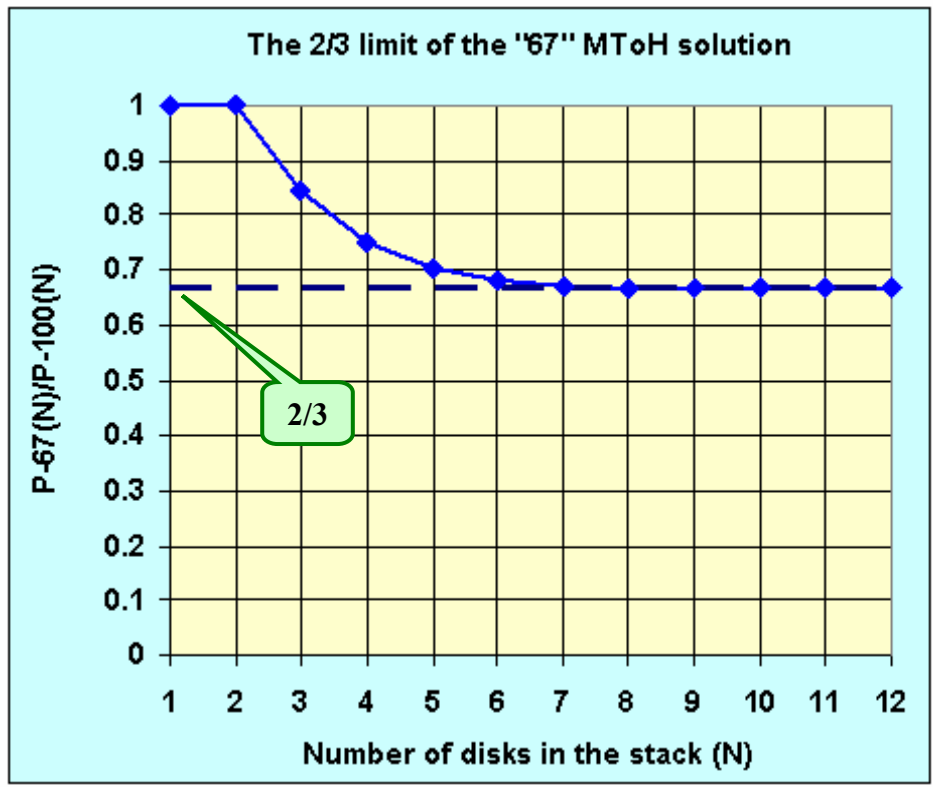

**Figure 8:** *"Efficiency" or "duration-ratio" curve for the "67" solution. As shown, the curve "quickly" approaches its limit of 2/3.*

All right. We have formulated a highly efficient solution, based essentially on the discovery that a three-disk **MToH** puzzle can be solved in just 11 moves. But did we find the most efficient solution? Is 2/3 the shortest relative-duration? Well, as was obvious right from the Abstract, the answer is "no". With a modified algorithm, triggered by new insights, the relative-duration limit can be pushed further down to 67/108 = ~ 62%.

### 2.4. The "62%" solution of the MToH puzzle

The "67" solution starts rather nicely. We efficiently move "down" N-1 disks to the Intermediate Post and move the N-th disk to its final rest on the Destination Post. But now, we either move a single disk or recursively move N-2 disks, using the $S_{100}$ Algorithm (see Table 5 and Equation 9). That is – on "folding" N-1 disks back up on the largest disk, we move the N-2 stack (four times) using the inefficient "100" Algorithm. As if the tower is permanently colored. As I will now show, a more efficient algorithm does exist.

As it turns out, on up-folding N-1 disks, we run into "SemiFree" States of the Tower. And a SemiFree Algorithm, to be discussed next, results in a shorter duration. Once we are done with the SemiFree Algorithm, we go back to the "62" Algorithm and swiftly complete it, enjoying what I think is the highest efficiency solution. Let's see then the definition of a SemiFree Tower and its associated disk-moving algorithm.

#### 2.4.1. The SemiFree Algorithm

On moving up N-1 disks (Over the largest disk) we run into a situation shown in Figure 9. The N-th disk is already on Post D which is now **Blue**, the N-1 disk is on Post S, which is "colored" **Red**, and we need to move N-3 disks onto Post S to clear the way for the N-2 disk to land **Red** on Post I. I discovered that moving the stack of N-3 disks from **Blue-D** to **Red-S** can be done rather efficiently. For example the reader can readily show that given the described Tower State, a stack of three disks (i.e. N-3 = 3) can be "relocated" in just 11 moves (and not 13). So we explore now this Tower State which we call "SemiFree".

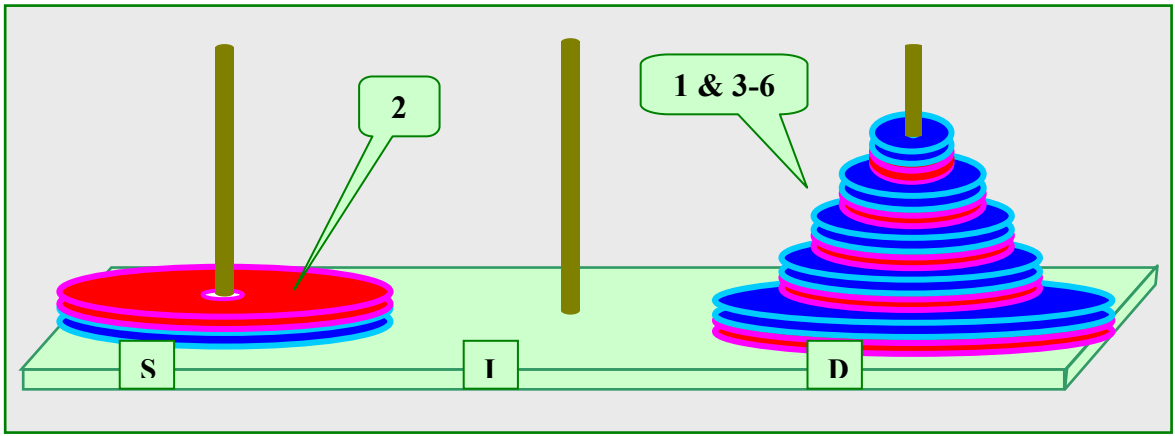

**Figure 9:** *An intermediate "SemiFree" State of the **MToH**. In the depicted example the top three disks on Post D are to be moved to Post S. This "mission" is efficiently accomplished by the SemiFree Algorithm.*

Our formal definition of a SemiFree Tower (consistent with the general definition in section 2.3.1) goes as follows:

An **MToH** is **SemiFree** if

- One of its posts – say – S, is permanently colored – say **Red** (by large disks)
- Another post – say – D, is permanently and **oppositely** colored (by large disks)
- The third post – I is **Free** so it has (at the start of the algorithm) a **Neutral** color (and can assume either color during execution of the algorithm)
- We need to move N disks from Post S to Post D using Post I

A SemiFree tower with N = 4 is shown in Figure 10.

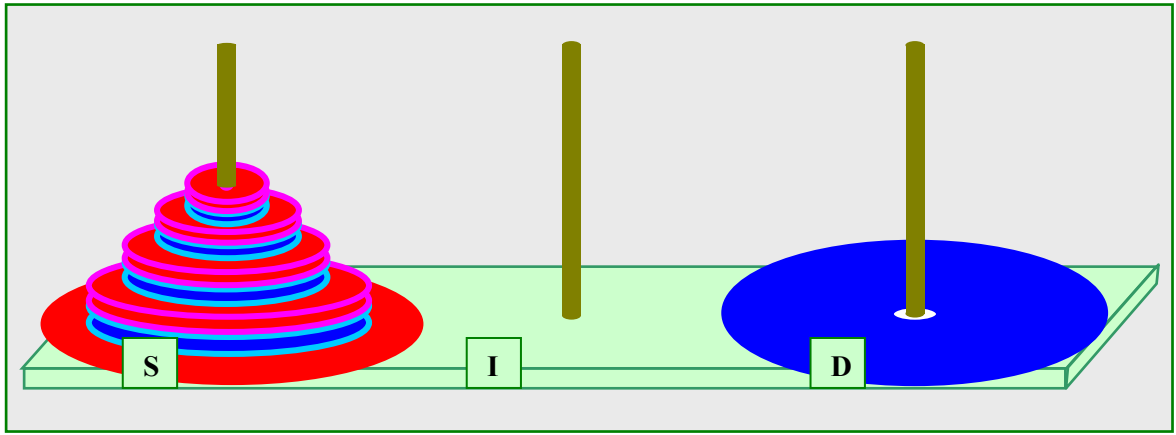

**Figure 10:** *Formal description of the SemiFree State of the MToH. Refer to the text for a rigorous definition. The mission here is to move the N disks now residing on Post S to reside on Post D. The mission is efficiently accomplished by the SemiFree Algorithm as described in the text.*

The SemiFree Algorithm is spelled-out by Table 7.

| Step # | Disks | From | To | # of moves | Comments |
|---|---|---|---|---|---|
| 1 | 3 to N | S (Red) | D (Blue) | $S_{SF}(N-2)$ | "SemiFree" **MToH** |
| 2 | 2 | S (Red) | I (Blue) | 1 | |
| 3 | 3 to N | D (Blue) | I (Blue) | $2*S_{100}(N-2)$ | |
| 4 | 1 | S (Red) | D (Blue) | 1 | |
| 5 | 3 to N | I (Blue) | D (Blue) | $2*S_{100}(N-2)$ | |
| 6 | 2 | I (Blue) | S (Red) | 1 | |
| 7 | 3 to N | D (Blue) | I (Red) | $1*S_{100}(N-2)$ | Post I changed color |
| 8 | 2 | S (Red) | D (Blue) | 1 | |
| 9 | 3 to N | I (Red) | D (Blue) | $1*S_{100}(N-2)$ | |

**Table 7:** *The SemiFree Algorithm. The algorithm moves N > 2 disks from S to D (through I), assuming the Source Post and the Destination Post are oppositely and permanently colored (in the actual solution both are occupied by larger disks).*

In terms of number of moves, we see from Table 7 -

$$S_{SF}(N) = S_{SF}(N-2) + 6 \cdot S_{100}(N-2) + \frac{3^2 - 1}{3^1 - 1}. \tag{14}$$

Equation 14 is recursive (the N-th value can be evaluated if the N-2 value is known). So it takes some effort to come up with closed-form expressions.

The closed-form expression for the number of moves of the k-th disk when executing the SemiFree Algorithm is given by Equation 15:

$$P_{SF}(k) = 2 \cdot 3^{k-2} + \frac{3^{k-1} - 3^2}{3^2 - 1} - \frac{3^{k-2} - 3^1}{3^2 - 1} + 1 \quad ; \text{ k odd}$$

$$P_{SF}(k) = 2 \cdot 3^{k-2} + \frac{3^{k-1} - 3^1}{3^2 - 1} - \frac{3^{k-2} - 3^2}{3^2 - 1} \quad ; \text{ k even} \tag{15}$$

The closed form expression for the total number of moves required to relocate a stack of N disks, executing the SemiFree Algorithm, is given by Equation 16:

$$S_{SF}(N) = (3^{N-1} + N - 1) + \frac{3^{N-1} - 3^0}{3^2 - 1} - \frac{N-1}{2} \quad ; \text{ N odd}$$

$$S_{SF}(N) = (3^{N-1} + N - 1) + \frac{3^{N-1} - 3^1}{3^2 - 1} - \frac{N-2}{2} \quad ; \text{ N even} \tag{16}$$

Table 8 lists the number of moves of the k-th disk for the first eight stack "heights".

As shown, the number of SemiFree moves is generally larger than the equivalent "67" number of moves (refer to Table 6) but is generally significantly smaller than the equivalent "100" number of moves (refer to Table 4).

| k \ N | 1 | 2 | 3 | 4 | 5 | 6 | 7 | 8 | SUM |
|---|---|---|---|---|---|---|---|---|---|
| 1 | 1 | | | | | | | | 1 |
| 2 | 1 | 3 | | | | | | | 4 |
| 3 | 1 | 3 | 7 | | | | | | 11 |
| 4 | 1 | 3 | 7 | 21 | | | | | 32 |
| 5 | 1 | 3 | 7 | 21 | 61 | | | | 93 |
| 6 | 1 | 3 | 7 | 21 | 61 | 183 | | | 276 |
| 7 | 1 | 3 | 7 | 21 | 61 | 183 | 547 | | 823 |
| 8 | 1 | 3 | 7 | 21 | 61 | 183 | 547 | 1641 | 2464 |

**Table 8:** *Number of disk-moves for the "SemiFree" Algorithm of the Magnetic Tower of Hanoi puzzle (Figure 10). N is the total number of disks to be moved from Post S to Post D, and k is the disk number in the ordered stack, counting from bottom to top. Equation 15 spells out the $P_{SF}(k)$ calculating expression and Equation 16 spells the expression for calculating the sum - $S_{SF}(N)$. Compare the numbers in this "SemiFree" table to the smaller corresponding numbers in the "67" table (Table 6) but Compare the numbers in this "SemiFree" table to the much larger corresponding numbers in the "100" table (Table 4), to realize the "SemiFree" savings.*

With Equation 16 in place, we can easily find the limit of the "duration-ratio" of "SF" vs. "100" as follows:

$$\frac{S_{SF}(N)}{S_{100}(N)} \xrightarrow{N \to \infty} \frac{3^{N-1} + 3^{N-1}/(3^2-1)}{3^N/2} = \frac{3}{4} \qquad (17)$$

The curve for the duration-ratio of "SF" vs. "100" (along with similar curves for "67" vs. "100" and "62" vs. "100") is shown in Figure 12 below.

Standing on top of the SemiFree "hill", we can already see the "62" summit. Let's then, following a short rest, climb the last mile.

### 2.4.2. The "62" Algorithm

With the SemiFree Algorithm in place, along with the "100" Algorithm and the "67" Algorithm, we now return to the original **MToH** and swiftly solve the puzzle.

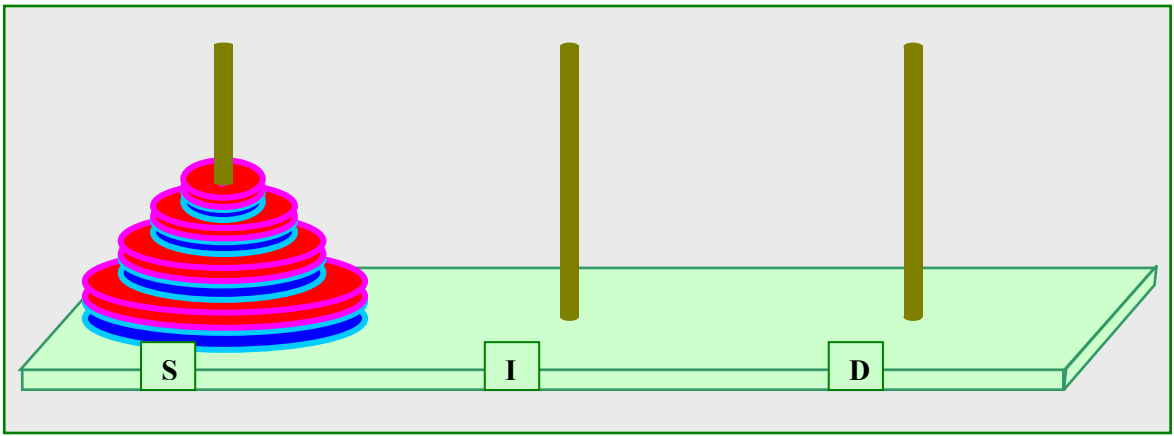

**Figure 11:** *The "regular" or "Free" **MToH** puzzle. The diagram shown in the figure is just a copy of the diagram shown in Figure 2, placed here for clarity and for reader's convenience.*

Figure 11 is just a copy of Figure 2, placed here for reader's convenience. Let' also repeat the game's objective - we want to efficiently relocate (i.e. relocate by a small number of moves) the N disks placed originally over the Source-Post onto the Destination-Post, subject to the Size Rule as well as to the Magnet Rule. To accomplish this mission (solve the puzzle efficiently), we present the "62" Algorithm.

Described in very general terms, the "62" Algorithm is made up of three steps –

- Move N-1 disks down onto Post I, colored **Blue** at the end of the sequence, using the "67" Algorithm

- Move two more disks to Post D while leaving N-3 disks on Post I, colored **Red** at the end of the sequence, using essentially the SemiFree Algorithm

- Move up the remaining N-3 disks (from Post I to Post D), using again the "67" Algorithm

An accurate, more detailed, description of the "62" Algorithm is given in Table 8.

| Step # | Disks | From | To | # of moves | Comments |
|---|---|---|---|---|---|
| 1 | 2 to N | S (Red) | I (Blue) | $S_{67}(N-1)$ | Going "down" |
| 2 | 1 | S (Red) | D (Blue) | 1 | |
| 3 | 3 to N | I (Blue) | D (Blue) | $2*S_{100}(N-2)$ | Start folding up here |
| 4 | 2 | I (Blue) | S (Red) | 1 | |
| 5 | 4 to N | D (Blue) | S (Red) | $S_{SF}(N-3)$ | SemiFree Algorithm |
| 6 | 3 | D (Blue) | I (Red) | 1 | Post I changed color |
| 7 | 4 to N | S (RED) | I (Red) | $2*S_{100}(N-3)$ | N-2 disks on Post I |
| 8 | 2 | S(RED) | D(Blue) | 1 | |
| 9 | 3 to N | I (Red) | D (Blue) | $S_{67}(N-2)$ | Efficient up-folding |

**Table 8:** *The "62" Algorithm for N ≥ 3. For N < 3, the "62" Algorithm coincides with the "67" Algorithm (see Table 5). The "62" Algorithm involves all three algorithms already analyzed – "100", "67", and "SemiFree".*
*Note that two "67" Algorithms are used in the "62" solution sequence. The one in step 1 is actually "67-Down" Algorithm. The one in step 9 is actually "67-Up" Algorithm. The "67-Up" Algorithm is a "time-reversed" "67-Down" Algorithm (and vice-versa – see Appendix 1). Necessarily, the move-counting equations (Equations 10 and 12) apply equally well to both algorithm variations.*

We want now to develop expressions for the number of puzzle-solving moves, for the "62" Algorithm. Looking at Table 8 we see only "recognized" algorithms ("100", "67", and "SemiFree"). Expression for the number of moves of the k-th disk for each of the three participating algorithms was already presented above. Similarly for the total number of moves. So now, for the "62" Algorithm, we simply sum the previously developed expressions -

$$P_{62}(k) = P_{67}(k) \; ; \; k \leq 3$$

$$P_{62}(k) = 2 \cdot P_{100}(k-2) + 2 \cdot P_{100}(k-3) +$$
$$+ P_{67}(k-1) + P_{67}(k-2) + P_{SF}(k-3) \; ; \; k > 3$$

(18)

And for the total number of moves -

$$S_{62}(N) = S_{67}(N) \; ; \; N \leq 3$$

$$S_{62}(N) = 2 \cdot S_{100}(N-2) + 2 \cdot S_{100}(N-3) +$$
$$+ S_{67}(N-1) + S_{67}(N-2) + S_{SF}(N-3) + \frac{3^2-1}{3^1-1} \; ; \; N > 3.$$

(19)

The two "67" Algorithms in Equation 19 are somewhat different. The first one, applied to N-1 disks, is actually "67-Down" Algorithm. The second one, applied to N-2 disks, is actually "67-Up" Algorithm. The "67-Up" Algorithm is a "time-reversed" "67-Down" Algorithm (and vice-versa – see Appendix 1). Necessarily, the move-counting equations related to the "67" Algorithm (Equations 10 and 12) apply equally well to both algorithms.

| N \ k | 1 | 2 | 3 | 4 | 5 | 6 | 7 | 8 | SUM |
|---|---|---|---|---|---|---|---|---|---|
| 1 | 1 | | | | | | | | 1 |
| 2 | 1 | 3 | | | | | | | 4 |
| 3 | 1 | 3 | 7 | | | | | | 11 |
| 4 | 1 | 3 | 7 | 19 | | | | | 30 |
| 5 | 1 | 3 | 7 | 19 | 53 | | | | 83 |
| 6 | 1 | 3 | 7 | 19 | 53 | 153 | | | 236 |
| 7 | 1 | 3 | 7 | 19 | 53 | 153 | 455 | | 691 |
| 8 | 1 | 3 | 7 | 19 | 53 | 153 | 455 | 1359 | 2050 |

**Table 9:** *Number of disk-moves for the "62" Algorith solving the Magnetic Tower of Hanoi puzzle (Figure 11, Equation 18 and Equation 19).*

Table 9 lists the number of moves of the k-th disk for the first eight stack "heights".

Looking at the number of moves for the "62" Algorithm as listed in Table 9, and comparing the numbers to the numbers listed in Table 6 for the "67" Algorithm, we indeed see some additional savings. For example, the total number of moves to solve the 8-disk **MToH** puzzle using the "67" Algorithm is 2194 while using the "62" Algorithm the number is only 2050. The "100" Algorithm, by the way, (the algorithm that solves a **C**olored-**MToH**), calls for 3280 moves (Table 4).

For the limit of the "duration-ratio" of "62" vs. "100", we retain the high N-powers of 3 to find –

$$\frac{S_{62}(N)}{S_{100}(N)} \xrightarrow{N \to \infty} \frac{2 \cdot 3^{N-2} + 2 \cdot 3^{N-3} + 3^{N-4} + \{3^{N-4}/(3^2 - 1)\}}{3^N/2} = \frac{67}{108} \quad (20)$$

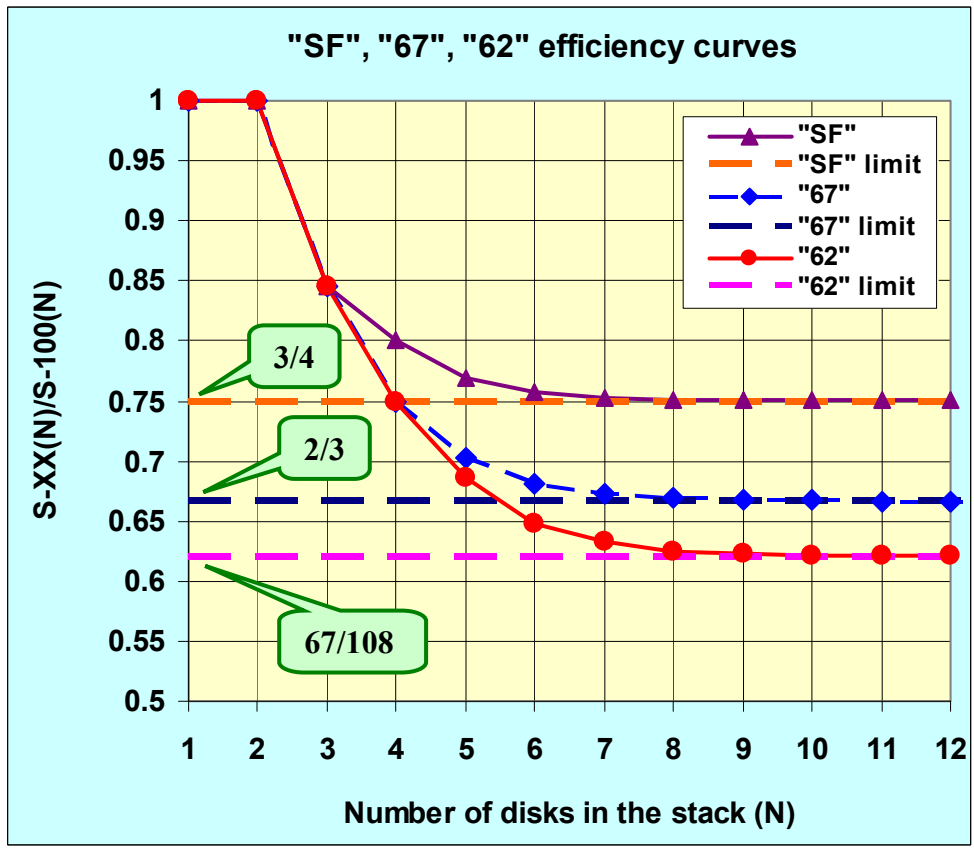

**Figure 12:** *Duration ratio curves for the "SF" Algorithm, the "67" Algorithm, and for the "62" Algorithm. The limits are 3/4, 2/3, and 67/108 respectively.*

We can see now the effect of "freedom" on solution "efficiency":

$$S_{SF}(N)/S_{100}(N) \xrightarrow{N \to \infty} 3/4 \qquad (21A)$$

$$S_{67}(N)/S_{100}(N) \xrightarrow{N \to \infty} 2/3 \qquad (21B)$$

$$S_{62}(N)/S_{100}(N) \xrightarrow{N \to \infty} 67/108 \qquad (21C)$$

The (fully) **Colored MToH** has the lowest efficiency defined as 1 [$P_{100}(N)/P_{100}(N) = 1$], the **SemiFree MToH** is next with efficiency of 3/4, and for the (completely) **Free MToH** we have found two algorithms – "67" and "62" with efficiencies of 2/3 and 67/108 respectively.

The three "duration-ratio" curves are shown in Figure 12.

Next is a set of recursive equations showing again the "starring" role of the number "3".

### 2.5. Recursive Relations

The four puzzle-solving algorithms discussed above follow remarkably simple (and thus elegant) recursive relations. This is true (recursive simplicity that is) for both the number of disk moves as well as for the total number of puzzle-solving moves. And in all relations, without exception, the leading term (on the right-hand side) is the product of the "current" value and the starring number - "3".

Let's see.

The "100" Algorithm:

$$P_{100}(k+1) = 3 \cdot P_{100}(k) \qquad (22A)$$

$$S_{100}(N+1) = 3 \cdot S_{100}(N) + 1 \qquad (22B)$$

The "67" Algorithm:

$$P_{67}(k+1) = 3 \cdot P_{67}(k) - 2 \quad ; \quad k \geq 2 \qquad (23A)$$

$$S_{67}(N+1) = 3 \cdot S_{67}(N) - 2 \cdot N + 3 \qquad (23B)$$

The "SF" Algorithm:

$$P_{SF}(k+1) = 3 \cdot P_{SF}(k) \quad ; \text{k odd} \qquad (24A)$$

$$P_{SF}(k+1) = 3 \cdot P_{SF}(k) - 2 \quad ; \text{k even} \qquad (24B)$$

$$S_{SF}(N+1) = 3 \cdot S_{SF}(N) - N + 2 \quad ; \text{N odd} \qquad (24C)$$

$$S_{SF}(N+1) = 3 \cdot S_{SF}(N) - N + 1 \quad ; \text{N even} \qquad (24D)$$

And the "62" Algorithm:

$$P_{62}(k+1) = 3 \cdot P_{62}(k) - 6 \quad ; \quad k \geq 4 \quad ; \text{k odd} \tag{25A}$$

$$P_{62}(k+1) = 3 \cdot P_{62}(k) - 4 \quad ; \quad k \geq 4 \quad ; \text{k even} \tag{25B}$$

$$S_{62}(N+1) = 3 \cdot P_{62}(N) - 5 \cdot (N-3) - 3 \quad ; \quad N \geq 3 \quad ; \text{N odd} \tag{25C}$$

$$S_{62}(N+1) = 3 \cdot P_{62}(N) - 5 \cdot (N-3) - 2 \quad ; \quad N \geq 3 \quad ; \text{N even} \tag{25D}$$

So much for the number of moves characterizing each of the four **MToH** puzzle-solving algorithms presented in this paper. Yet, before concluding, I wish to bring forward another section. A "Color-Crossings" section that is. The section presents the color of each of the three posts during the entire solving procedure, in a **graphical form**. Shedding colorful light onto the **MToH** puzzle.

Let's see.

### 2.6. Color-Crossings

To visualize color-crossings, I asked the computer to record the color of each post for each move, from start to finish, and designate each color by a number - "1" for Red, "0" for Neutral and "-1" for Blue.

Selected recordings are shown by the two figures below.

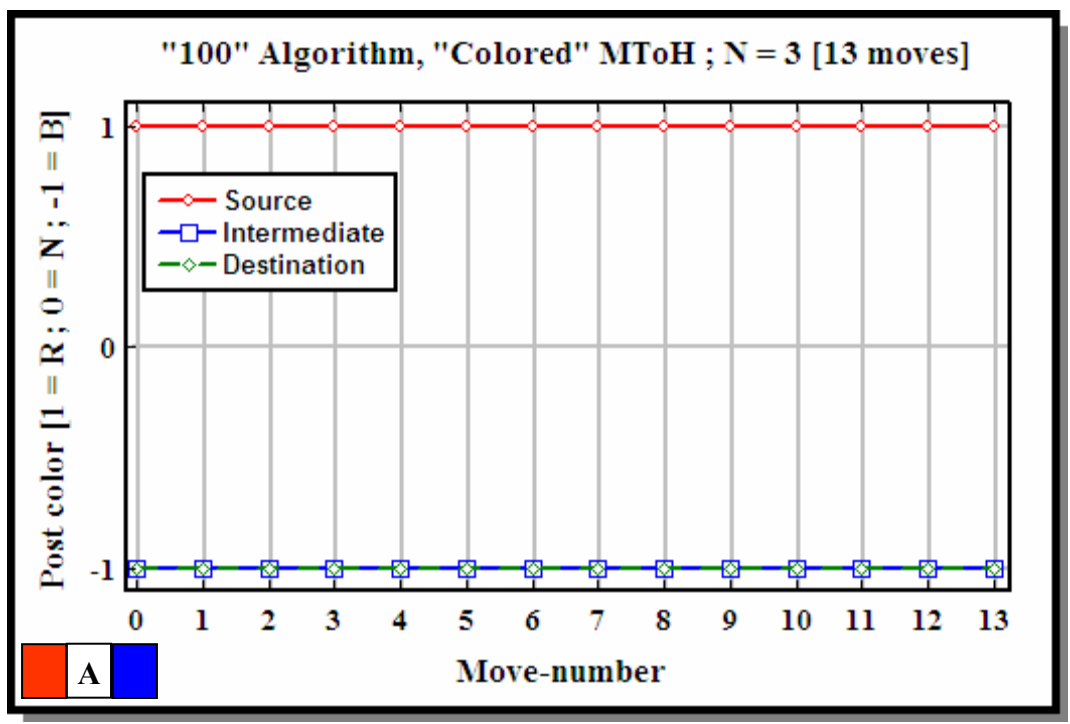

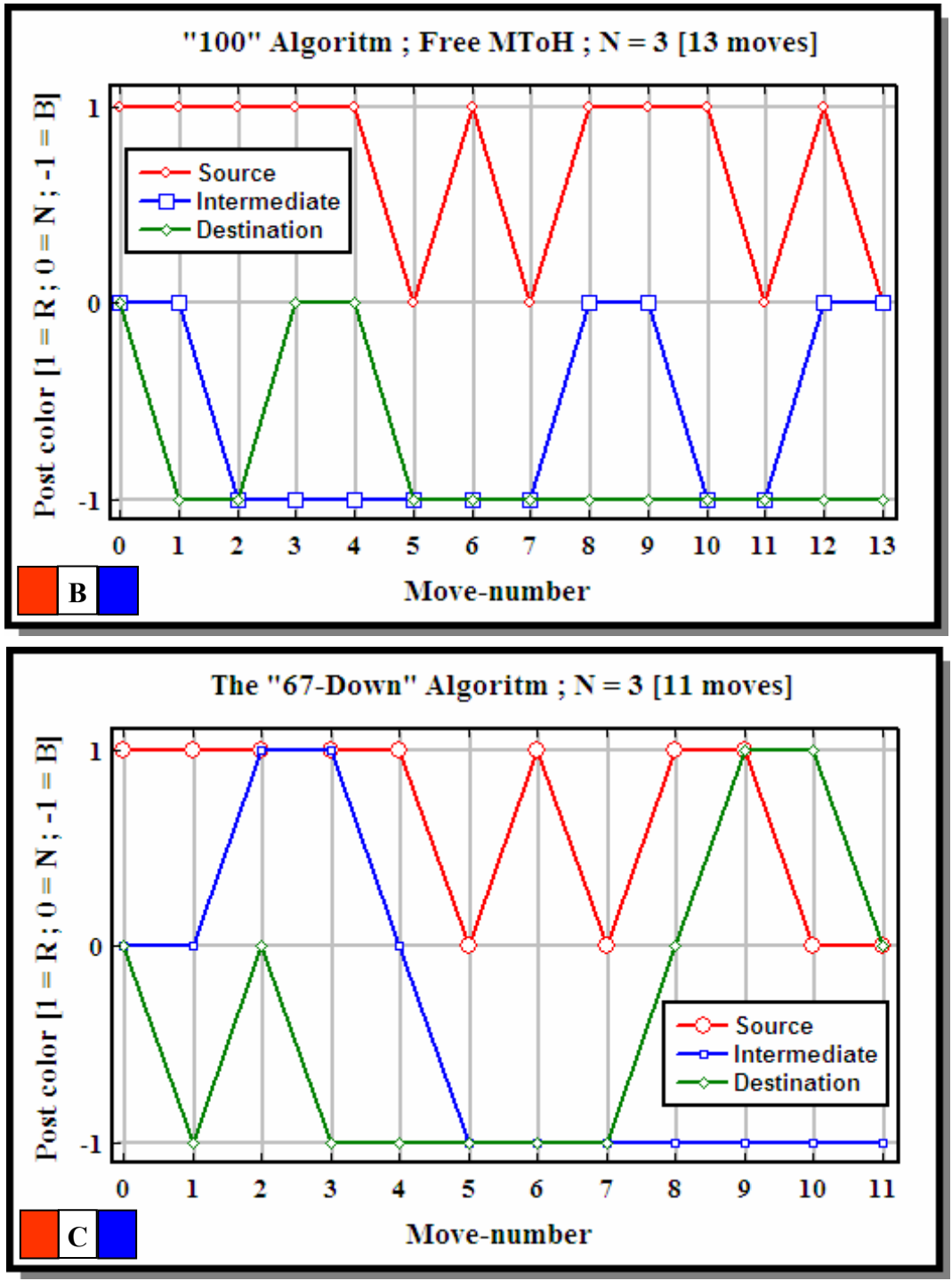

**Figure 13:** *Color-Crossings charts. All three charts are associated with a height 3 **MToH**. A – Colored **MToH** and the "100" Algorithm. 13 moves, No Color-Crossing. B – Free **MToH** and the "100" Algorithm. Still 13 moves. Still no Color-Crossing (see text for details). C - Free **MToH** and the "67-Down" Algorithm. Two Color-Crossings, only 11 moves.*

The three charts of Figure 13 all relate to height 3 of the **MToH**.

The top one (A) shows the color of each of the posts for a **C**olored **MToH** (**CMToH**). In this case of a "Permanently Colored" Tower, the posts are pre-colored **Red**-**Blue**-**Blue** and the "100" Algorithm curves, not surprising, stay horizontal throughout the entire 13-move solution.

The middle one (B) relates to a "regular" or "Free" **MToH**, solved by exactly the same "100" Algoritm as was the case for 13A. Now, during the 13-move solution, we see each of the three posts wonders between Neutral and one color, **never crossing Neutral to "visit" the opposite color**.

The bottom one (C) relates to the "62-Down" Algorithm, solving the **MToH** puzzle in, as we know very well by now, just 11 moves. In this case we see two Color-Crossings. By "Color-Crossing" I refer to a move sequence where a post goes from one color through Neutral (and may stay there for a short while) to the **opposite** color. Such Color-Crossing is exersiced by the Intermediade Post in moves 3,4 and 5 and by the Destination Post in moves 7,8 and 9 of the "62-Down" Algorithm. **These Color-Crossings "take responsibility" for the shorter-duration solution** of only 11 moves.

Next, Figure 14 shows Color-Crossing charts for an **MToH** of height 5, comparing the crossings of the "67-Down" Algorithm (A) to the crossings of the "62" Algoritm (B).

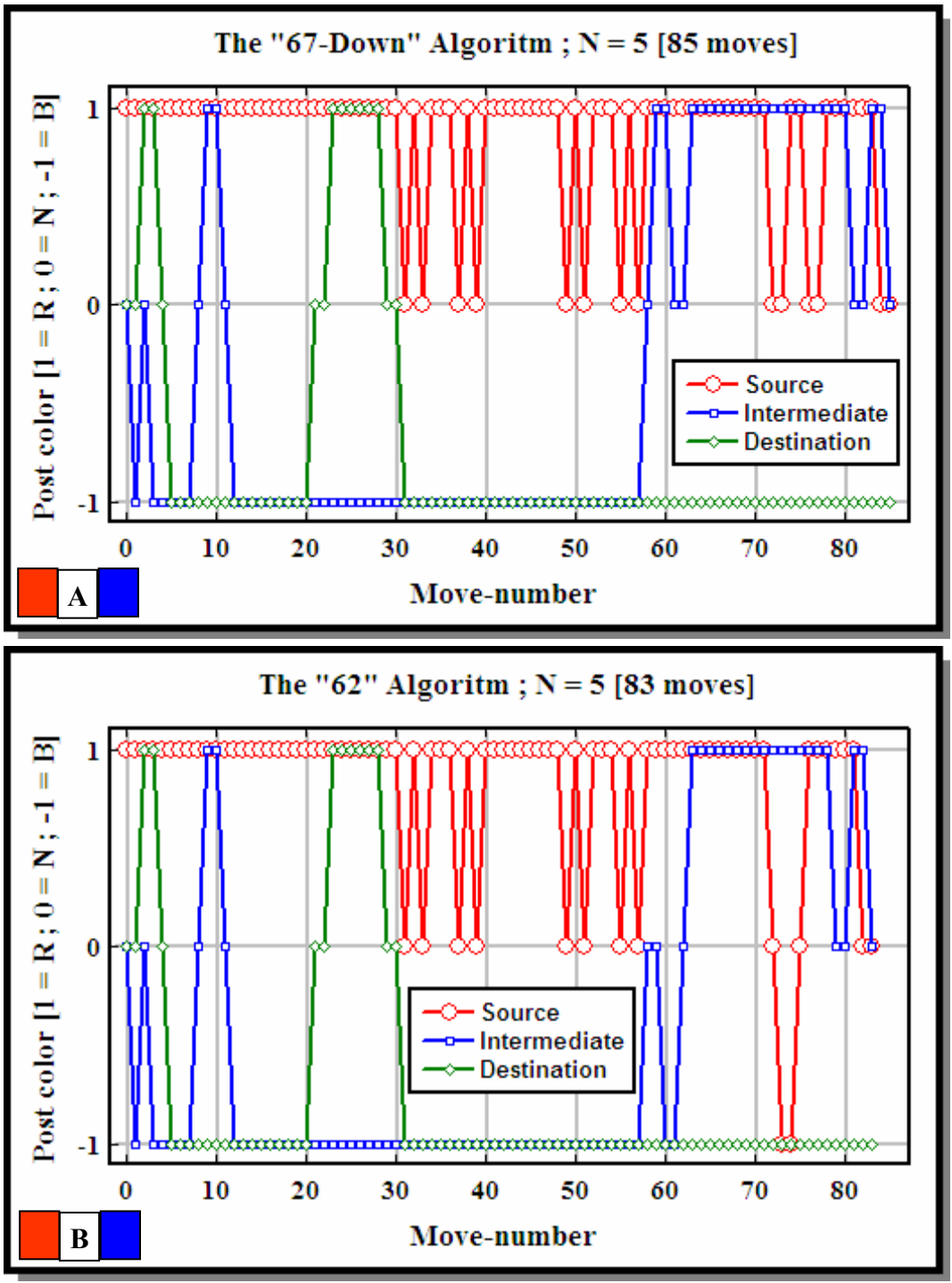

**Figure 14:** *Color-Crossings charts for a Free **MToH** of height 5. A – "67-Down" Algorithm. Six Color-Crossings (see Table 10 below). 85 moves. B – "62" Algorithm. Eight Color-Crossings (see Table 10 below). Eighty-three moves.*

Comparing the top chart (Figure 5A) to the bottom chart (Figure 5B), We see additional two Color-Crossings of the Source-Post (71 through 73 ; 74 through 76) for the "62" Algorithm (Figure 5B). Again we wittness the correlation between larger number of Color-Crossings and a solution of a shorter duration.

To see this Crossings-Duration correlation, we listed in Table 10 the number of Color-Crossings of each post for the first eight **MToH** heights.

| N | 1 | 2 | 3 | 4 | 5 | 6 | 7 | 8 |
|---|---|---|---|---|---|---|---|---|
| **67-Down-S** | 0 | 0 | 0 | 0 | 0 | 0 | 0 | 0 |
| **67-Down-I** | 0 | 0 | 1 | 2 | 3 | 4 | 5 | 6 |
| **67-Down-D** | 0 | 0 | 1 | 2 | 3 | 4 | 5 | 6 |
| **67-Down-total** | 0 | 0 | 2 | 4 | 6 | 8 | 10 | 12 |
| P67-Down(N) | 1 | 4 | 11 | 30 | 85 | 248 | 735 | 2194 |
|  |  |  |  |  |  |  |  |  |
| **67-Up-S** | 0 | 0 | 0 | 2 | 2 | 4 | 4 | 6 |
| **67-Up-I** | 0 | 0 | 1 | 1 | 3 | 3 | 5 | 5 |
| **67-Up-D** | 0 | 0 | 0 | 0 | 0 | 0 | 0 | 0 |
| **67-Up-total** | 0 | 0 | 1 | 3 | 5 | 7 | 9 | 11 |
| P67-Up(N) | 1 | 4 | 11 | 30 | 85 | 248 | 735 | 2194 |
|  |  |  |  |  |  |  |  |  |
| **62-S** | 0 | 0 | 0 | 1 | 2 | 2 | 4 | 4 |
| **62-I** | 0 | 0 | 1 | 2 | 3 | 8 | 9 | 14 |
| **62-D** | 0 | 0 | 0 | 2 | 3 | 4 | 5 | 6 |
| **62-total** | 0 | 0 | 1 | 5 | 8 | 14 | 18 | 24 |
| P62(N) | 1 | 4 | 11 | 30 | 83 | 236 | 691 | 2050 |

**Table 10:** *Color-Crossings for three algorithms for the first eight **MToH** heights. For "high" Towers, the posts in the "62" Algorithm make significantly larger numbers of Color-Crossings vs. the corresponding numbers for the two "67" Algorithms.*

We did that for three Algorithms – "67-Down", "67-Up" and "62". We had to split the "67" Algorithm because, as shown, the Color-Crossing pattern for the "67-Down" Algorithm differs slightly from the Color-Crossing pattern for the "67-Up" Algorithm. Both, however, solve the **MToH** puzzle in exactly the same number of moves. And while both are characterized, for each stack height, by a similar number of crossings, they both display significantly smaller number of Color-Crossings (for "high" stacks) when compared to the number of Color-Crossimgs of the "62" Algorithm. And we know that the "62" Algorithm solution is of shorter duration. For high stacks then, the correlation discovered and discussed in relation to height 3, and height 5 persists.

So much for the **MToH** move analysis.

Now just a few organizing remarks before concluding.

All four Algorithms discussed above – "100", "67", "SemiFree" and "62", are recursive. Explicit recursive functions that run on NUMERIT[5] ("Mathematical & Scientific Computing") are listed in Appendix 1. Also listed in Appendix 1 are "program managing" functions that were written for program clarity and for better program managability.

A "movie" showing the "62" Algorithm solving a height five **MToH** in (**only**) 83 moves can be seen here[6].

A few pictures, showing actual realization of the **MToH** puzzle, are shown in Appendix 2.

Let's conclude now.

### 3. Concluding remarks

The task of the "Monks of Hanoi" is nearing completion. **The big disk has been moved.** Evidently, $2^{63}$ = 9.223372036854775808*$10^{18}$ seconds have already past since the Monks started performing their routine (always without the slightest hesitation). So **SOON** *"the world will end"* [1]! If only the *command of the ancient prophecy* would have been to move the 64 disks under the rules of the **Magnetic Tower of Hanoi**. If that was the case, we would still have
$((3^{64}-1)/2)*(67/108) – 2^{63}$ = 3.550259505549357568*$10^{29}$ seconds of **color**ful life ahead of us (out of the original 1.06507785166480704*$10^{30}$ seconds since they started). But let's not worry. Let's enjoy our world, with the innovations it offers, for the remaining 9.223372036854775808*$10^{18}$ seconds.

"Always without the slightest hesitation". I used this phrase in the previous paragraph. Because as a matter of fact, for the classical base-2 **ToH,** determinism prevails. If the play moves "forward" (on the down-sequence for example, N-2 disks go over the freshly moved disk number N-1 and not back over disk number N) then **the moves are mandatory**. No need to think, no reason to hesitate. The same applies to the **C**olored **M**agnetic **T**ower **o**f **H**anoi. True, both Towers span their respective bases perfectly, but the puzzle solution has an element of monotony in it. Solving (efficiently) the Free Magnetic Tower of Hanoi puzzle is a different story. On one hand, when counting moves, the number "3" stars. If you look back through this paper, you will find this number (3) in all of the equations from Equation 3 and on. Without exception. In some early equations implicitely. These early "hints" do not decieve us. As we easily realize now - "1" is actually $3^0$ ; "4" is actually $(3^2 - 3^0)/(3^1 - 3^0)$ ; "11" is actually $3^{(3-1)} + (3 - 1)$. And so, indeed, number "3" is everywhere. However, not only number 3 stars, but the game is intricate. The puzzle solution may progess in more than one path. The puzzle presents more than one option to the player. The Tower therefore calls for thinking, justifies hesitation. It is **Freedom** that makes the **MToH** puzzle so **colorful**.

## 4. References

[1]  http://en.wikipedia.org/wiki/%C3%89douard_Lucas

[2]  http://en.wikipedia.org/wiki/Tower_of_Hanoi

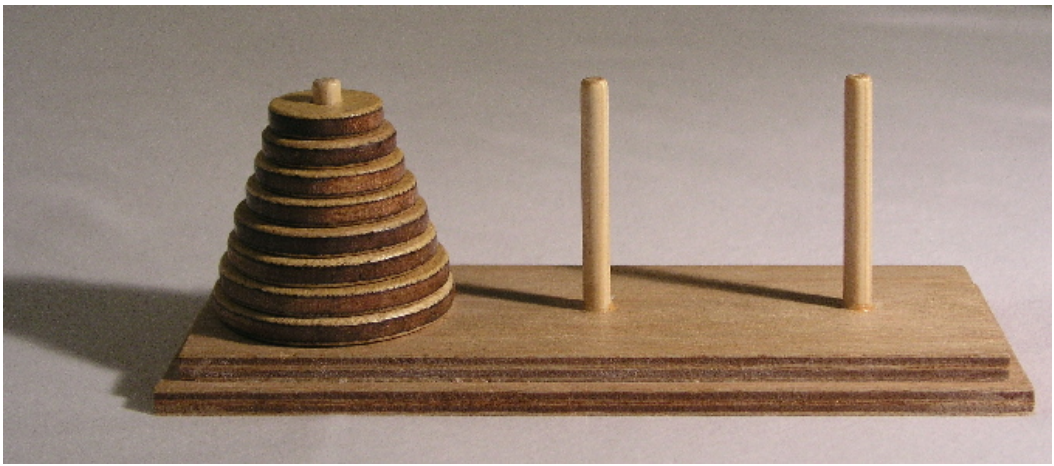

A model set of the Towers of Hanoi (with 8 disks)

*...If the legend were true, and if the priests were able to move disks at a rate of one per second, using the smallest number of moves, it would take them $2^{64} - 1$ seconds or roughly 600 billion years[1]; it would take 18,446,744,073,709,551,615 turns to finish.*

[3]  http://www.cut-the-knot.org/recurrence/hanoi.shtml

[4]  "The Magnetic Tower of Hanoi", Uri Levy, Journal of Recreational Mathematics 35:3, to be published (~April 2010)

[5]  http://www.numerit.com/

[6]    *A "movie" showing the "62" Algorithm solving a height five **MToH** in (only) 83 moves:*
       http://www.numerit.com/maghanoi/
       User name: maghanoi
       Password: mtoh2009

**Appendix 1:** Recursive functions for the "62" solution

Listed in this Appendix are all the functions by which the "62" Algorithm solves the **MToH** puzzle. The functions run on NUMERIT[4] – a "Mathematical & Scientific Computing" environment. Five of the functions (D,G,H,I,J) are recursive (call themselves). These functions, the "heart" of the game, may offer important clues needed to decipher the **MToH** puzzle.

**A.**
**solve_MToH_puzzle(n,s,d,i)**

**B.**
**function solve_MToH_puzzle(n,s,d,i)**
    if n=1
        move_down_67(n,s,d,i)
    return
    if n=2
        move_down_67(n,s,d,i)
    return
    move_down_67(n-1,s,i,d)
    move(n,s,d)
    move_all_but_n_up(n,i,d,s)

**C.**
**function move_all_but_n_up(n,i,d,s)**
    if n > 2
        move_busy_2and1(n-2,i,s,d)
        move_busy_1and2(n-2,s,d,i)
        move(n-1,i,s)
        move_semifree_BNR(n-3,d,s,i)
        move(n-2,d,i)
        move_busy_2and1(n-3,s,d,i)
        move_busy_1and2(n-3,d,i,s)
        move(n-1,s,d)
        move_up_67(n-2,i,d,s)

**D.**
**function move_semifree_BNR(n,s,d,i)**
    if n > 2
        move_semifree_BNR(n-2,s,d,i)
        move(n-1,s,i)
        move_busy_2and1(n-2,d,s,i)
        move_busy_1and2(n-2,s,i,d)
        move(n,s,d)
        move_busy_2and1(n-2,i,s,d)
        move_busy_1and2(n-2,s,d,i)
        move(n-1,i,s)
        move_busy_1and2(n-2,d,i,s)
        move(n-1,s,d)
        move_busy_2and1(n-2,i,d,s)
    return
    if n = 1
        move_busy_2and1(1,s,d,i)
    return
    if n = 2
        move_busy_2and1(2,s,d,i)
    return

**E.**
**function move_down_67_3disks(n,s,d,i)**
    if n > 0
        move_busy_2and1(n-1,s,i,d)
        move(n,s,d)
        move_busy_1and2(n-2,i,s,d)
        move_busy_2and1(n-2,s,d,i)
        move(n-1,i,s)
        move_busy_2and1(n-2,d,i,s)
        move(n-1,s,d)
        move_busy_1and2(n-2,i,d,s)

**F.**
**function** move_up_67_3disks(n,s,d,i)
    if n > 0
        move_busy_2and1(n-2,s,i,d)
        move(n-1,s,d)
        move_busy_1and2(n-2,i,s,d)
        move(n-1,d,i)
        move_busy_1and2(n-2,s,d,i)
        move_busy_2and1(n-2,d,i,s)
        move(n,s,d)
        move_busy_2and1(n-1,i,d,s)

**G.**
**function** move_down_67(n,s,d,i)
    if n > 3
        move_down_67(n-1,s,i,d)
        move(n,s,d)
        move_busy_2and1(n-2,i,s,d)
        move_busy_1and2(n-2,s,d,i)
        move(n-1,i,s)
        move_busy_1and2(n-2,d,i,s)
        move(n-1,s,d)
        move_busy_2and1(n-2,i,d,s)
    return
    if n=2
        move_busy_1and2(2,s,d,i)
    return
    if n=1
        move_busy_1and2(1,s,d,i)
    return
    move_down_67_3disks(3,s,d,i)

**H.**
**function move_up_67(n,s,d,i)**
    if n > 3
        move_busy_1and2(n-2,s,i,d)
        move(n-1,s,d)
        move_busy_2and1(n-2,i,s,d)
        move(n-1,d,i)
        move_busy_2and1(n-2,s,d,i)
        move_busy_1and2(n-2,d,i,s)
        move(n,s,d)
        move_up_67(n-1,i,d,s)
    return
    if n=2
        move_busy_1and2(2,s,d,i)
    return
    if n=1
        move_busy_1and2(1,s,d,i)
    return
    move_up_67_3disks(3,s,d,i)

**I.**
**function move_busy_2and1(n,i,s,d)**
    if n > 0
        move_busy_2and1(n-1,i,s,d)
        move_busy_1and2(n-1,s,d,i)
        move(n,i,s)
        move_busy_2and1(n-1,d,s,i)

**J.**
**function move_busy_1and2(n,s,d,i)**
    if n > 0
        move_busy_1and2(n-1,s,i,d)
        move(n,s,d)
        move_busy_2and1(n-1,i,s,d)
        move_busy_1and2(n-1,s,d,i)

**Appendix 2:** Realization of the Magnetic Tower of Hanoi

Below – a few pictures of a realized version of the Magnetic Tower of Hanoi.

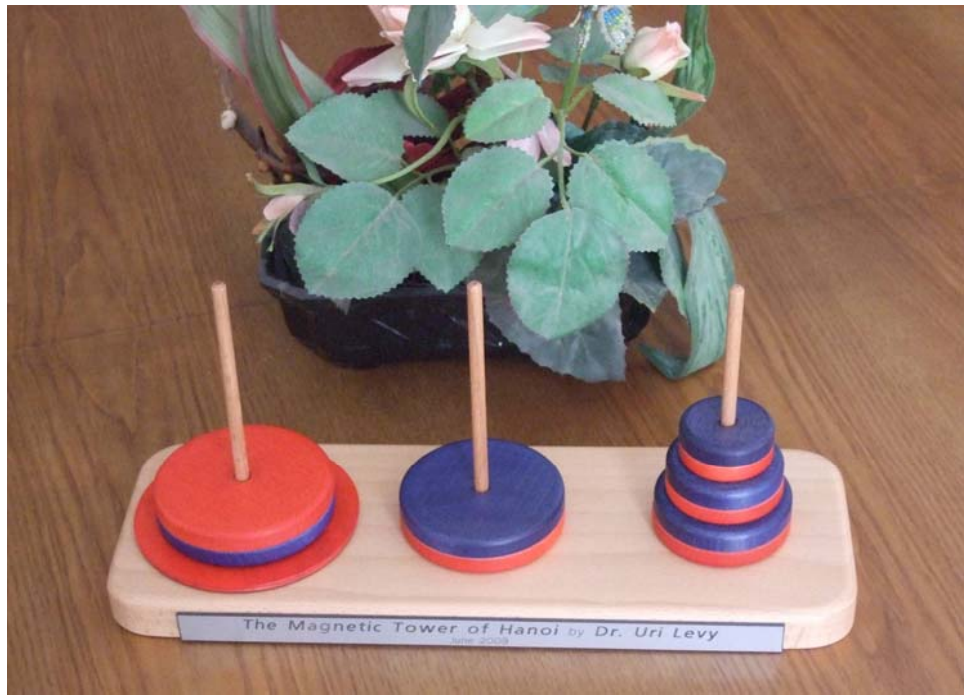

Figure A2-1: *MToH puzzle realized.*

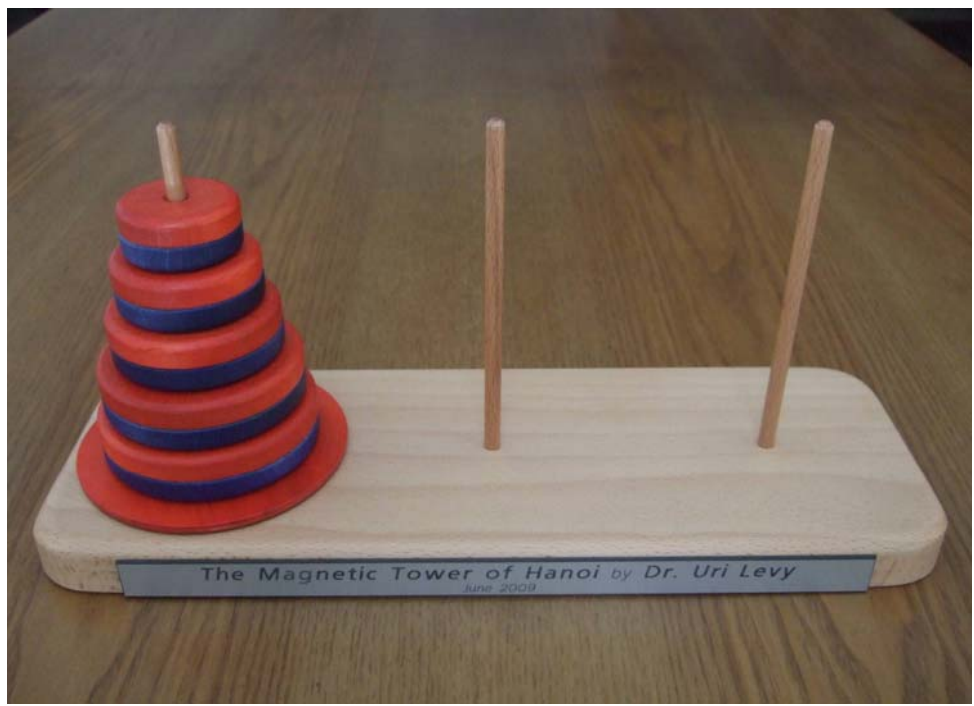

Figure A2-2: *The "free" or "classical" MToH.*

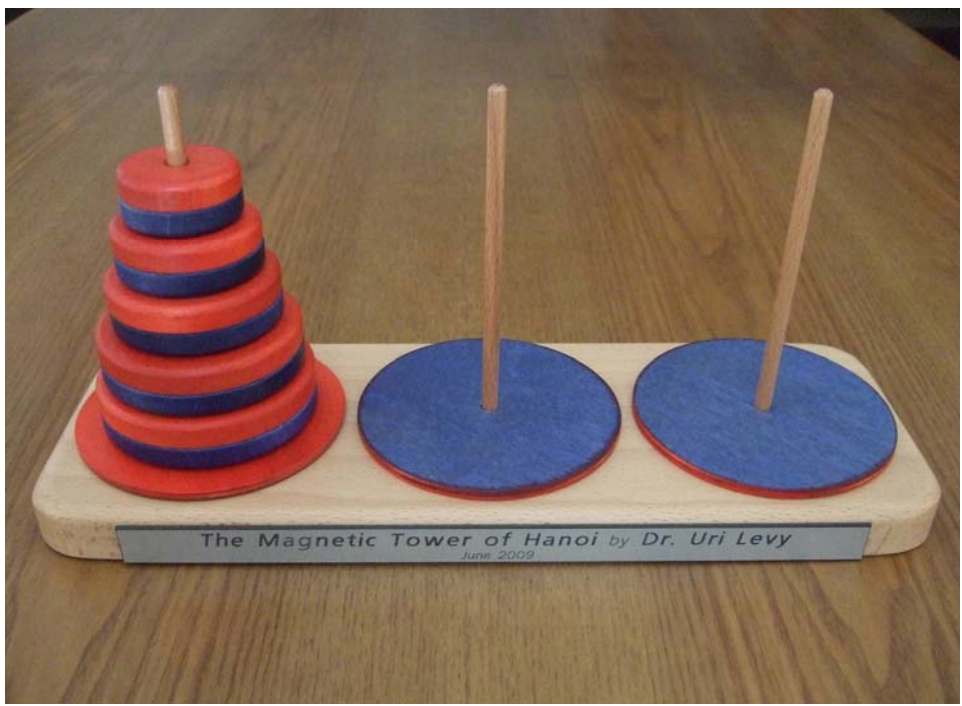

Figure A2-3: The Colored *MToH* puzzle realized.

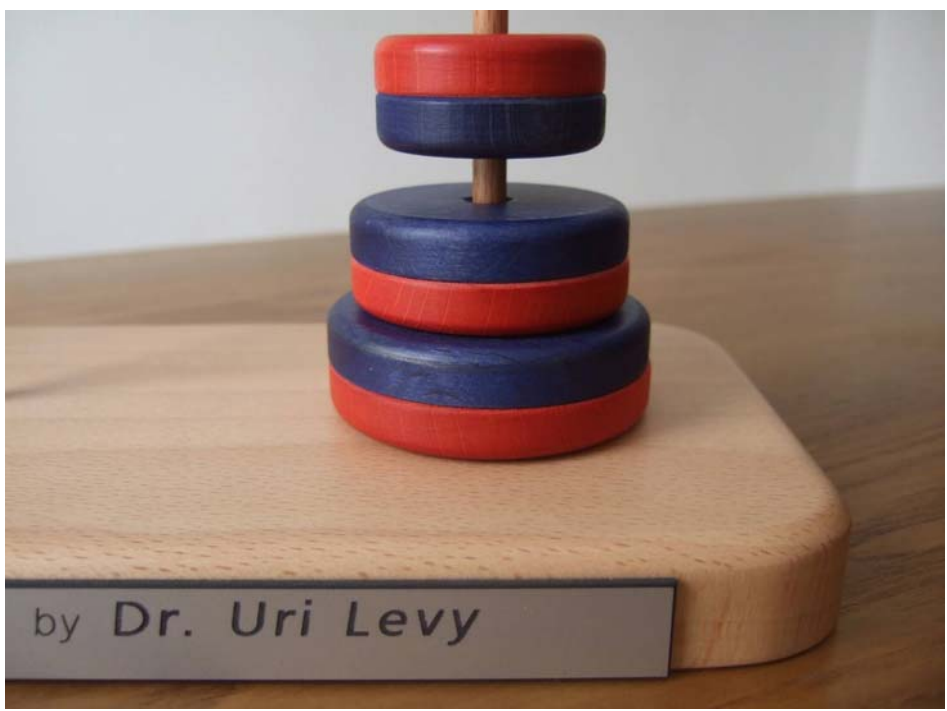

Figure A2-4: Illegal move. *The hidden inside-disk magnets physically prevent complete disk placement, indicating illegal move.*